\documentstyle[12pt]{article}
\begin{document}
\newcommand{\Si}{\Sigma}
\newcommand{\tr}{{\rm tr}}
\newcommand{\ad}{{\rm ad}}
\newcommand{\Ad}{{\rm Ad}}
\newcommand{\ti}[1]{\tilde{#1}}
\newcommand{\om}{\omega}
\newcommand{\Om}{\Omega}
\newcommand{\de}{\delta}
\newcommand{\al}{\alpha}
\newcommand{\te}{\theta}
\newcommand{\vth}{\vartheta}
\newcommand{\be}{\beta}
\newcommand{\lm}{\lambda}
\newcommand{\La}{\Lambda}
\newcommand{\D}{\Delta}
\newcommand{\ve}{\varepsilon}
\newcommand{\vf}{\varphi}
\newcommand{\G}{\Gamma}
\newcommand{\ip}{\hat{\upsilon}}
\newcommand{\Ip}{\hat{\Upsilon}}
\newcommand{\ga}{\gamma}
\newcommand{\ze}{\zeta}
\newcommand{\li}{\lim_{n\rightarrow \infty}}
\newcommand{\mat}[4]{\left(\begin{array}{rr}{#1}&{#2}\\{#3}&{#4}
\end{array}\right)}
\newcommand{\si}{\sigma}
\newcommand{\beq}[1]{\begin{equation}\label{#1}}
\newcommand{\eq}{\end{equation}}
\newcommand{\beqn}[1]{\begin{eqnarray}\label{#1}}
\newcommand{\eqn}{\end{eqnarray}}
\newcommand{\p}{\partial}
\newcommand{\di}{{\rm diag}}
\newcommand{\su}{{\bf su_2}}
\newcommand{\uo}{{\bf u_1}}
\newcommand{\ep}{\epsilon}
\newcommand{\GL}{{\rm GL}(N,{\bf C})}
\newcommand{\gl}{gl(N,{\bf C})}
\newcommand{\rar}{\rightarrow}
\newcommand{\upar}{\uparrow}
\newcommand{\sm}{\setminus}
\newcommand{\ms}{\mapsto}
\newcommand{\bp}{\bar{\partial}}
\newcommand{\sect}[1]{\setcounter{equation}{0}\section{#1}}
\renewcommand{\theequation}{\thesection.\arabic{equation}}
\newtheorem{predl}{Proposition}[section]
\newtheorem{defi}{Definition}[section]
\newtheorem{rem}{Remark}[section]
\newtheorem{cor}{Corollary}[section]
\newtheorem{lem}{Lemma}[section]
\newtheorem{theor}{Theorem}[section]
\def\smd{\mathbin{>\mkern -8mu\lhd}}

\begin{titlepage}
\setcounter{footnote}0
\begin{center}
\hfill ITEP TH-7/96\\
\hfill MPI-96-60
\hfill hep-th/9605005\\
\vspace{0.3in}

{\LARGE\bf Double coset construction of moduli space of holomorphic
bundles and  Hitchin systems.}\\
\vspace{0.15in}

\bigskip

\bigskip

{\Large A.Levin\footnote{E-mail address: alevin@mpim-bonn.mpg.de
}}$\phantom{hj}^{\dag}$,
{\Large M.Olshanetsky\footnote{E-mail address:
olshanez@vxdesy.desy.de}}$\phantom{hj}^{\ddag}$,
\\
\bigskip

{\sf Max-Planck-Institut f\"{u}r Mathematik, Bonn}\\
\bigskip

\bigskip

\begin{quotation}{
$\phantom{hj}^{\dag}$ --
{\it On leave from International Institute for Nonlinear Studies at 
Landau Inst, Vorob'iovskoe sch. 2, Moscow, 117940, Russia}

$\phantom{hj}^{\ddag}$ --
{\it On leave from ITEP, Bol.Cheremushkinskaya, 25, 
Moscow, 117 259, Russia}\\
}
\end{quotation}
\end{center}
\end{titlepage}

\title{}
\author{}
\maketitle
\date{}



\begin{abstract}
We present a description of the moduli space of holomorphic vector
bundles over Riemann curves as a double coset space which is differ
from the standard loop group construction. Our approach is based on
equivalent definitions of holomorphic bundles, based on the 
transition maps
or on the first order differential operators. Using this approach we
present  two independent derivations of the Hitchin integrable systems.
We define a "superfree" upstairs systems from which Hitchin systems are
obtained by  three step hamiltonian reductions. A special attention is 
being given on the Schottky parameterization of curves. 
 \end{abstract}

\section{Introduction}
The moduli space of holomorphic vector bundles over Riemann surfaces are
popular subject in algebraic geometry and  number theory. In
mathematical physics they were investigated due to relations with the
Yang-Mills theory \cite{AB} and the Wess-Zumino-Witten theory
\cite{KZ,B}. The conformal blocks in the  WZW theory satisfy the
Ward identities which take a form of differential equations on the
moduli space \cite{TUY,BSch}. In this approach the moduli space is
 described
as a double coset space of a loop group defined on a small circle 
on a Riemann surface \cite{PS}. 

The main goal of the paper is an alternative description of the moduli
space and the Hitchin integrable systems \cite{H} based on this construction.
 We start with a special group valued
field on a Riemann surface which is defined as a map from a holomorphic 
basis in a
vector bundle to a $C^{\infty}$ basis. 
This field is an analogous of the tetrade field in 
the General Relativity and 
we call it the Generalized Tetrade Field (GTF). 
The holomorphic structures 
can be extracted from GTF. They are
described via the holomorphic transition maps, or by means of
 the  operators $d''$ . 
The former are invariant under the action of 
the global $C^{\infty}$ transformations, while the later 
under the action 
of the local holomorphic transformations.
   It  allows to define the moduli space as a double
coset space of GTF with respect to the actions of the local holomorphic
 transformations and the global $C^{\infty}$ transformations. 

 We introduce a cotangent
bundle to GTF and invariant symplectic structure on it. The 
cotangent bundle to the moduli of holomorphic bundles can be obtained by
the symplectic 
factorizations over the action of two types of 
commuting gauge transformations. 
 This cotangent bundle is a phase space of the Hitchin integrable 
 systems \cite{H}. 
The tetrade fields in their turns are sections of the principle bundle
 over the Riemann
surface, which satisfy some constraints equations. We interpret them as
 moment constraints in a big "superfree system"
 with a special gauge symmetry. 
This space is a cotangent bundle to the principle bundle.
Thus the Hitchin systems are obtained by the three step 
symplectic reductions from this space. 

We investigate specially our reductions in terms
 of Schottky parameterization,
which is a particular case of the general construction. 
This parameterization was 
used to derived the Knizhnik-Zamolodchikov-Bernard equations on the higher
genus curves \cite{B,Lo,I}. On the other hand the
quantum second order Hitchin Hamiltonians coincide with them on the
critical level. 

\section{Moduli of holomorphic vector bundles}
\setcounter{equation}{0}

Let $\Si=\Si_{g}$ be a nondegenerate Riemann curve of genus $g$ with
$g>1$ .
We will consider in this section a set of stable holomorphic structures  on 
complex vector bundles  
over $\Si$ \cite{AB}. To define them  we proceed in two ways 
based on the \^{C}ech and the Dolbeault cohomologies. 
Eventually, we come to the moduli space ${\cal L}$ of stable holomorphic
bundles over $\Si_{g}$ and represent them as a double coset space 
(Proposition 2.3).
\bigskip

{\bf 1.} Consider  a vector bundle $V$ over $\Si_{g}$ .
To be more concrete we assume 
that the 
structure group of $V$ is $GL(N,{\bf C})$.
Let ${\cal U}_a,~a=1,\ldots$ be a covering of $\Si_{g}$ by open subsets.
We consider two bases in
$V$ the holomorphic $\{e^{hol}\}$ basis and  the smooth
$C^{\infty}$ $\{e^{C^{\infty}}\}$ one.
 In local coordinates 
$(z_a\in{\cal U}_a)$ $$e^{hol}_a=e^{hol}(z_a),~
e^{C^{\infty}}_a=e^{C^{\infty}}(z_a,\bar{z}_a).$$
Let $h$ be the transition map between them
$h_a=h(z_a,\bar{z}_a)$.
Then locally in ${\cal U}_a$ we have
\beq{a1}
h_ae^{C^{\infty}}_a=e_a^{hol}.
\eq
We can consider $h_a$ as the sections $\Om_{C^{\infty}}^0({\cal U}_a , P)$ 
of the adjoint bundle
$P=$Aut $V$. We  call the field $h$  a 
generalized tetrade field (GTF).
It follows from the definitions of the bases that there exists a global
section for $ e_{C^{\infty}}$
\beq{a3}
 e^{C^{\infty}}_a(z_a,\bar{z}_a)= e^{C^{\infty}}_b(z_b(z_a),
\bar{z}_b(\bar{z}_a)),~
z_a\in{\cal U}_{ab}={\cal U}_a\cap{\cal U}_b\neq\emptyset,
\eq
where $z_b=z_b(z_a)$ are holomorphic functions 
defining a complex structure
on $\Si_g$. On the other hand the
 transformations of $e^{hol}$ are holomorphic maps
\beq{a4}
e^{hol}_a(z_a)=g_{ab}(z_a)e^{hol}_b(z_b(z_a)),~~g_{ba}(z_b)=g_{ab}^{-1}(z_a(z_b))
\eq
$$
g_{ab}
\in \Om_{hol}^0({\cal U}_{ab},{\rm Aut}~ V), 
~~(\bp g_{ab}=0,~\bp=\p_{\bar z_a}).
$$
These matrix functions define the holomorphic structure
 in the vector bundle $V$.

 We can describe the same holomorphic structure working with the smooth
basis $e^{C^{\infty}}$ in $V$. Let  
\beq{a5}
\bar{A}_a=h_a^{-1}\bar{\p}h_a.
\eq
Then the basis $ e^{C^{\infty}}$ is annihilated by the operator 
$d''_A|_{{\cal U}_{a}}=\bp+\bar A_a$
$$
(\bp+\bar A_a)e^{C^{\infty}}_a=0.
$$

The GTF transformations $h$ in (\ref{a1}) by no means free.
Let ${\cal R}_\Si$ be the subset of sections in $P$
 which satisfies the following
conditions
\beq{a6}
{\cal R}_\Si=\{h\in\Om_{C^{\infty}}^0({\cal U}_a, P)~|
~h_a^{-1}\bar{\p}h_a|_{{\cal U}_{ab}}=
h_b^{-1}\bar{\p}h_b|_{{\cal U}_{ab}}, ~\forall~
 {\cal U}_{ab}\neq\emptyset,~a,b=1,\ldots\},
\eq
$$
(\bar{A}_a(z_a)=\bar{A}_b(z_b(z_a)),~z_a\in{\cal U}_{ab}).
$$
\bigskip
\begin{predl}
Conditions (\ref{a1}) and (\ref{a6}) are equivalent.
\end{predl}
{\it Proof.}
Since $e^{C^{\infty}}_b=e^{C^{\infty}}_a$ in ${\cal U}_{ab}$
 (\ref{a1}) implies
\beq{a7}
g_{ab}=h_ah_b^{-1}.
\eq
Then the holomorphicity of $g_{ab}$ implies (\ref{a6}) .
If $h\in{\cal R}_\Si$, then
(\ref{a7}) defines the transition map for some holomorphic basis $e^{hol}$.
The basis $e^{hol}h$ satisfies (\ref{a3}) and therefore can be taken as
$e^{C^{\infty}}$.
$\Box$
\bigskip

Consider the group
\beq{a8}
{\cal G}_{\Si}=\{\Om_{C^{\infty}}^0({\cal U}_a,P),~a=1,\ldots\}.
\eq
It transforms local basses of $V$ over ${\cal U}_a$.
The group acts on itself by the left 
and right multiplications.

There are two   subgroups of ${\cal G}_{\Si}$. Let 
$x_a\in\Om_{C^{\infty}}^0({\cal U}_a,P)$. Then
\beq{a10}
{\cal G}_{\Si}^{hol}=
\{x_a\rar f_ax_a~|~f\in\Om^0_{hol}(\Si,P)\},
\eq
\beq{a11}
{\cal G}_{\Si}^{C^{\infty}}=\{x_a\rar x_a\vf_a~|~
\vf\in \Om_{C^{\infty}}^0(\Si,P),~ 
\vf(z_b(z_a),\bar{z}_b(\bar{z}_a))= \vf(z_a,\bar{z}_a)~
z_a\in{\cal U}_{ab}\}.
\eq
We can consider the GTF (\ref{a6}) as a subset in ${\cal G}_{\Si}$.
We have the following evident statement
\bigskip 
\begin{predl}
The left and right actions of ${\cal G}_{\Si}^{hol}$ and 
${\cal G}_{\Si}^{C^{\infty}}$
 leave invariant  ${\cal R}_{\Si}$. 
\end{predl}
\bigskip
In other words
$$
\def\normalbaselines{\baselineskip20pt
       \lineskip3pt    \lineskiplimit3pt}
\def\mapright#1{\smash{
        \mathop{\longrightarrow}\limits^{#1}}}
\def\mapup#1{\Big\uparrow\rlap
       {$\vcenter{\hbox{$\scriptstyle#1$}}$}}
\matrix{
{\cal G}_{\Si}:~&{\cal G}_{\Si}&\mapright{right~ mltpl.}&{\cal G}_{\Si}\cr
\cup             &\cup          &                       &\cup    \cr
{\cal G}_{\Si}^{C^{\infty}}:~
&{\cal R}_{\Si}&\mapright{right ~mltpl.}&{\cal R}_{\Si}\cr}
$$
and
$$
\def\normalbaselines{\baselineskip20pt
       \lineskip3pt    \lineskiplimit3pt}
\def\mapright#1{\smash{
        \mathop{\longrightarrow}\limits^{#1}}}
\def\mapup#1{\Big\uparrow\rlap
       {$\vcenter{\hbox{$\scriptstyle#1$}}$}}
\matrix{
{\cal G}_{\Si}:~&{\cal G}_{\Si}&\mapright{left~ mltpl.}&{\cal G}_{\Si}\cr
\cup             &\cup          &                       &\cup    \cr
{\cal G}_{\Si}^{hol}:~&{\cal R}_{\Si}&\mapright{left~ mltpl.}&{\cal R}_{\Si}\cr}
$$

{\bf 2}.
Consider the space of holomorphic structures on the  bundles $V$ and $P$.
Since \\ $g>1$  there is an open subset of stable holomorphic
structures. The  holomorphic structures can be defined in two  ways.
In the first type  of the construction,  which we call the D-type,   
the  holomorphic structures are defined by
the covariant operators. For $V$ they are
 $$
d''_A:\Om_{C^{\infty}}^{0}(\Si,V)\rar\Om_{C^{\infty}}^{0,1}(\Si,V).
$$
It means that $\bar{A}$ satisfies (\ref{a6}).
The  holomorphic structure is consistent with the complex structure on $\Si_{g}$
 such that for any section
$s\in\Om_{C^{\infty}}^{0}(\Si,V)$ and $f\in C^{\infty}(\Si)$ 
$d''_A(fs)=(\bar{\p}f)s+fd''_As.$
The space of holomorphic structures ${\cal L}_{\Si}^D$ on $P$ 
is defined in the similar way 
\beq{a12}
{\cal L}_{\Si}^D=
\{d''_A=\bar{\p}+\bar{A}:\Om_{C^{\infty}} ^0(\Si,P)\rar
\Om_{C^{\infty}}^{(0,1)}(\Si,P)\}.
\eq
with the action in the adjoint representation. The stable holomorphic structures 
${\cal L}_{\Si}^{D,st}$
are an open subset
in (\ref{a12}). The automorphisms of the holomorphic structures are given
by the action of the gauge group ${\cal G}_{\Si}^{C^{\infty}}$ (\ref{a11})
\beq{a12a}
d''_A\rar \vf^{-1}d''_A\vf,~\vf\in{\cal G}_{\Si}^{C^{\infty}}.
\eq
They preserve the subset ${\cal L}_{\Si}^{D,st}$. 
{\it The moduli space ${\cal L}$ of stable holomorphic structures} on $P$ 
is the quotient space
\beq{a13}
{\cal L}={\cal L}_{\Si}^{D,st}/{\cal G}_{\Si}^{C^{\infty}}.
\eq
It is a smooth complex manifold with tangent space at $\bar{A}$ 
is isomorphic to \\
$H^{(0,1)}(\Si,{\rm Lie}(\GL))$.
Its dimension is given by the Riemann-Roch theorem
\beq{a14}
\dim{\cal L}=N^2(g-1)+1.
\eq
The left action of the gauge transformations ${\cal G}_{\Si}^{hol}$ (\ref{a10})
 does not change \\
$\bar{A}_a=h_a^{-1}\bp h_a,~a=1,\ldots$. Therefore the space 
${\cal L}_{\Si}^{D}$ (\ref{a12}) can be represented as the quotient space
${\cal L}_{\Si}^{D}={\cal G}_{\Si}^{hol}\backslash {\cal R}_{\Si}$.
There is an open subset in ${\cal R}_{\Si}^{st}$ such that
the subset of the stable holomorphic structures is the quotient space
$$
{\cal L}_{\Si}^{D,st}={\cal G}_{\Si}^{hol}\backslash{\cal R}_{\Si}^{st}.
$$
The main statement of this section follows immediately from (\ref{a13})
\bigskip
\begin{predl}
The moduli space ${\cal L}$ of stable holomorphic structures on $P$ can be
represented as the double  
 coset space
\beq{a15}
{\cal L}=
{\cal G}_{\Si}^{hol}\backslash {\cal R}_{\Si}^{st}/{\cal G}_{\Si}^{C^{\infty}}.
\eq
\end{predl}
\bigskip 

{\bf 3.} An alternative  description of the holomorphic structures 
 in terms of the \^{C}ech  cohomologies, which  we  call 
 the C-type construction is
based on the transition maps (\ref{a4}), (\ref{a7}).
The collection of  transition maps 
\beq{a16}
{\cal L}_{\Si}^{Ch}=\{g_{ab}(z_a)=h_a(z_a)h^{-1}_b(z_b(z_a)),~z_a\in
{\cal U}_{ab},~a,b=1,\ldots ,\}.
\eq
defines the holomorphic structures on $V$ or $P$ depending on the choice
of the representations. Again we choose the open subset of stable holomorphic
structures ${\cal L}_{\Si}^{C,st}$ in ${\cal L}_{\Si}^{Ch}$.
The gauge group ${\cal G}_{\Si}^{hol}$ acts as the automorphisms of
 ${\cal L}_{\Si}^{C,st}$
\beq{a17}
g_{ab}\rar f_ag_{ab}f^{-1}_b,~f_a=f(z_a),~
f_b=f_b(z_b(z_a)),~
f\in{\cal G}_{\Si}^{hol}.
\eq

The space ${\cal L}_{\Si}^{Ch}$ has a transparent description in terms of graphs.
Consider the skeleton of the covering $\{{\cal U}_a,~a=1,\ldots\}$. It is an
oriented graph,
whose vertices are some fixed inner points in ${\cal U}_a$ and edges $L_{ab}$ 
connect those $V_a$ and $V_b$ for which $U_{ab}\neq\emptyset$. 
We choose an orientation of the graph, saying that $a>b$ on the edge
 $L_{ab}$ and put
the holomorphic function $z_b(z_a)$ which defines the holomorphic map
from ${\cal U}_a$ to ${\cal U}_b$.  
Then the space ${\cal L}_{\Si}^{Ch}$ can
be defined by the following data. To each edge $L_{ab},~a>b$ we attach a matrix
valued function $g_{ab}\in\GL$ along with $z_b(z_a)$. The gauge fields 
$f_a$ are living on the 
vertices $V_a$ and the gauge transformation is (\ref{a17}).

The moduli space of stable holomorphic bundles 
is defined as the factor space under this action
\beq{a18}
{\cal L}={\cal G}_{\Si}^{hol}\backslash{\cal L}_{\Si}^{Ch,st}.
\eq
The tangent space to the moduli space in this approach is 
$H^1(\Si,{\rm Lie}(\GL))$ extracted from the \^{C}ech complex.
Though ${\cal L}_{\Si}^{Ch,st}$ differs from ${\cal L}_{\Si}^{D,st}$ we obtain
 the same moduli space ${\cal L}$ of stable holomorphic structures on $P$ due 
to the equivalence of the Dolbeault and the \^{C}ech cohomologies.

In this construction the right action of  ${\cal G}_{\Si}^{C^{\infty}}$ 
(\ref{a11}) leaves the transition maps $g_{ab}$ invariant. Therefore 
\beq{a19}
{\cal L}_{\ti{\Si}}^{Ch,st}={\cal R}_{\Si}^{st}/{\cal G}_{\Si}^{C^{\infty}}.
\eq
Taking into account (\ref{a18}) we come to the same construction of the
moduli space as the double coset space (\ref{a15}).
\bigskip

{\bf 4.} We fit the components of our construction in the 
 exact bicomplex
\begin{center}
${\bf M}$
\end{center}
$$
\def\normalbaselines{\baselineskip20pt
       \lineskip3pt    \lineskiplimit3pt}
\def\mapright#1{\smash{
        \mathop{\longrightarrow}\limits^{#1}}}
\def\mapup#1{\Big\uparrow\rlap
       {$\vcenter{\hbox{$\scriptstyle#1$}}$}}
\matrix{
& &\mapup{\p^C}&   &\mapup{\p^C}&   &\mapup{\p^C}&\cr
0&\mapright{}&\Om_{hol}^0({\cal U}_{ab}, P)&\mapright{i} 
&\Om_{C^{\infty}}^0({\cal U}_{ab},P) & \mapright{} & 
\Om_{C^{\infty}}^{(0,1)}({\cal U}_{ab},{\rm End}V)& \mapright{\bp}\cr
& &\mapup{\p^C}&   &\mapup{\p^C}&   &\mapup{\p^C}&\cr
0&\mapright{}&\Om_{hol}^0({\cal U}_{a},P)
&\mapright{i} 
&\Om_{C^{\infty}}^0({\cal U}_{a},P) & \mapright{} & 
\Om_{C^{\infty}}^{(0,1)}({\cal U}_{a},{\rm End}V)
& \mapright{\bp}\cr
& &\mapup{i}&   &\mapup{i}&   &\mapup{i}&\cr
0&\mapright{}&\Om_{hol}^0(\Si,P)&\mapright{i} 
&\Om_{C^{\infty}}^0(\Si,P) & \mapright{} & 
\Om_{C^{\infty}}^{(0,1)}(\Si,{\rm End}V)\}&
 \mapright{\bp}\cr
& &\mapup{}&   &\mapup{}&   &\mapup{}&\cr
& &0& &0& &0&\cr}
$$

Here $\p^C$ are the \^{C}ech differentials, $i$ are 
the augmentations. The right arrows from
 $\Om_{C^{\infty}}^0(*,P)$ to $\Om_{C^{\infty}}^{(0,1)}(*,{\rm End}V)$
are of the type $h\rar h^{-1}\bp h$. We have
$$
g_{ab}\in\Om_{C^{\infty}}^{(0,1)}({\cal U}_{ab},{\rm End}V),~
h_a\in\Om_{C^{\infty}}^0({\cal U}_{a},P),
$$
$$
\de\bar{A}\in\Om_{C^{\infty}}^{(0,1)}({\cal U}_{a},{\rm End}V).
$$
If these fields satisfy the tetrade conditions (\ref{a1}),(\ref{a4}) or
(\ref{a6}) then they lie in the images of $i$. The Dolbeault cohomologies 
$H^{(0,1)}(\Si,{\rm End}V)$ that define the tangent space to the moduli
space are living in $\Om_{C^{\infty}}^{(0,1)}(\Si,{\rm End}V)$ and the
\^{C}ech cohomologies $H^1(\Si,{\rm End}V)$ 
in $\Om_{hol}^0({\cal U}_{ab},{\rm End}V)$. 
Their equivalence can be derived from the properties of the double spectral
sequence.

The gauge transformations also can be incorporated in the  exact bicomplex

\begin{center}
${\bf G}$:
\end{center}
$$
\def\normalbaselines{\baselineskip20pt
       \lineskip3pt    \lineskiplimit3pt}
\def\mapright#1{\smash{
        \mathop{\longrightarrow}\limits^{#1}}}
\def\mapup#1{\Big\uparrow\rlap
       {$\vcenter{\hbox{$\scriptstyle#1$}}$}}
\matrix{
&&\mapup{\p^C}&   &\mapup{\p^C}&   &\mapup{\p^C}&\cr
0&\mapright{}&\Om_{hol}^0({\cal U}_{ab},{\rm End}P)&\mapright{i} 
&\Om_{C^{\infty}}^0({\cal U}_{ab},{\rm End}P) & \mapright{\bp} & 
\Om_{C^{\infty}}^{(0,1)}({\cal U}_{ab},{\rm End}P)
& \mapright{\bp}\cr
&&\mapup{\p^C}&   &\mapup{\p^C}&   &\mapup{\p^C}&\cr
0&\mapright{}&\Om_{hol}^0({\cal U}_{a},{\rm End}P)&\mapright{i} 
&\Om_{C^{\infty}}^0({\cal U}_{a},{\rm End}P) & \mapright{\bp} & 
\Om_{C^{\infty}}^{(0,1)}({\cal U}_{a},{\rm End}P)
& \mapright{\bp}\cr
&&\mapup{i}&   &\mapup{i}&   &\mapup{i}&\cr
0&\mapright{}&\Om_{hol}^0(\Si,{\rm End}P)&\mapright{i} 
&\Om_{C^{\infty}}^0(\Si,{\rm End}P) & \mapright{\bp} & 
\Om_{C^{\infty}}^{(0,1)}(\Si,{\rm End}P)
& \mapright{\bp}\cr
& &\mapup{}&   &\mapup{}&   &\mapup{}&\cr
& &0& &0& &0&\cr}
$$

Let 
$\ep^{hol}\in{\rm Lie}({\cal G}_{\Si}^{hol}),
~\ep^{C^{\infty}}\in{\rm Lie}({\cal G}_{\Si}^{C^{\infty}}).$ Then
$$
\ep^{hol}\in{\rm Image~of} (\Om_{hol}^0({\cal U}_{a},{\rm End}P))~
{\rm in}~\Om_{C^{\infty}}^0({\cal U}_{a},{\rm End}P),
$$  
$$
\ep^{C^{\infty}}\in{\rm Image~of} (\Om_{C^{\infty}}^0(\Si,{\rm End}P)~
{\rm in}~
\Om_{C^{\infty}}^0({\cal U}_{a},{\rm End}P).
$$
 The actions of the gauge group
(see (\ref{a12a}) and (\ref{a17}))
\beq{a20}
\de^{C^{\infty}}\bar{A}_a=\bp\ep^{C^{\infty}}_a+[\bar{A}_a,\ep^{C^{\infty}}_a],
\eq
\beq{a21}
\de^{hol}g_{ab}=\ep^{hol}_ag_{ab}-g_{ab}\ep^{hol}_b.
\eq
More generally,  ${\bf M}$ is the bigraded  ${\bf G}$ module.
The action of ${\bf G}$ is consistent with  the both differential $\p^C$
and $\bp$. The differentiations  take into account the bigradings of 
${\bf M}$ and ${\bf G}$.
 The actions (\ref{a20}),(\ref{a21}) are particular cases of these actions .

\section{ The Schottky specialization.}
\setcounter{equation}{0}

We apply the general scheme to the particular covering of $\Si_{g}$
based on the Schottky parameterization. 
Consider the Riemann sphere with $2g$ circles
 ${\cal A}_a,{\cal A}'_a,~
a=1,\ldots g$. Each circle lies in the 
external part of others. Let $\ga_a$ be $g$  projective
maps ${\cal A'}_a=\ga_a {\cal A}_a,~\ga_a\in PSL(2)$. The Schottky group 
$\G$ is a free group
 generated by $\ga_a,a=1\ldots g$.
The exterior part of all the circles 
$$
\ti{\Si}={\bf P}^1/\cup_{b=1}^{2g}D_b
$$ is the fundamental domain of $\G$.
The surface $\Si$ is obtained from $\ti{\Si}$  by the pairwise gluing
of the circles ${\cal A'}_a=\ga_a {\cal A}_a$ and
 ${\cal A}_a$.
We have only one nonsimpliconnected 2d cell ${\cal U}_a\sim\ti{\Si}$
with selfintersections ${\cal U}_{aa'}={\rm vicinity}~{\cal A}_a
={\rm vicinity}~{\cal A'}_a$. 
We choose  $g$  local coordinates $z_a,~a=1,\ldots,g$, which define 
the parameterizations of the  internal disks of ${\cal A}_a$ circles.
In this case the holomorphic maps can be written as
$z_{a'}(z_a)=\ga_a(z_a)$. The GTF ${\cal R}_\Si$ (\ref{a6}) is a twisted field $h$
on  $\ti{\Si}$
$$~h^{-1}\bar{\p}h(z_a,\bar{z}_a)=
h^{-1}\bar{\p}h(\ga_a(z_a),\overline{\ga_a(z_a)}), ~a=1,\ldots,g.
$$ 
In the definition of ${\cal G}_{\Si}^{C^{\infty}}$ (\ref{a11}) "the
periodicity conditions" take the form
$$
\vf(\ga_a(z_a),\bar{\ga}_a(\bar{z}_a))= \vf(z_a,\bar{z}_a),~
z_a\in{\rm vicinity~of~}{\cal A}_{a}.
$$
The transition maps (\ref{a4}),(\ref{a7}) defining ${\cal L}_{\Si}^{Ch}$ are 
\beq{a22}
g_a=g_{aa'}(z_a)=
h(z_a,\bar{z}_a)h^{-1}(\ga_a(z_a),\bar{\ga}_a(\bar{z}_a)),~
a=1,\ldots,g.
\eq
The gauge group ${\cal G}_{\Si}^{hol}$ acts as a global holomorphic 
transformation on $\ti{\Si}$. In the local coordinates we have
\beq{a23}
g_{a}\rar f_ag_{a}f^{-1}_{\ga_a},~f_a=f(z_a),~g_a=g_a(z_a),~
f_{\ga_a}=f(\ga_a(z_a)).
\eq
  In local coordinates $g_a$ have the  form of Laurent polynomials.
$g_a(z_a)=\sum g_a^{(k)}z_a^k.$
Thus  in
 this parameterization the set of holomorphic structures on the vector bundles 
${\cal L}_{\Si}^{Ch}$ 
 can be 
identified with  the collection of the loop groups
$L(GL_a(N,{\bf C}))$. But in fact, taking into account the adjoint action of
the gauge group
(\ref{a23}),  one concludes that the precise form of components
 is the semidirect product 
$ L(\GL)\smd PSL(2)=\{g(z)\smd\ga(z)\}.$  Thus
\beq{a24}
{\cal L}_{\Si}^{Ch}=\oplus_{a=1}^g   L_a(\GL)\smd PSL(2)_a,
\eq
where the subgroups $\{PSL(2)_a\},a=1\ldots g$ are responsible for 
the complex structure on $\Si$.
To define the stable bundles one should choose an open subset in 
$L_a(GL(N,{\bf C}))$. 

Consider the bundles over genus $g=1$ curves. Though the bundles are unstable 
this case  can be completely described in the wellknown 
terms. The Schottky parameterization means the realization
of elliptic curve as an annulus. Let
$~\ga(z)=qz,~q=\exp(2\pi i\tau )$. The holomorphic bundles are defined by 
  the loop group extended by the shift operator
\beq{a25}
{\cal L}_{\Si}^{Ch}= L(GL(N,{\bf C}))\smd\exp(2\pi i\tau z\p)  .
\eq
 The gauge action (\ref{a23}) 
$$g(z)\rar f(z)g(z)f^{-1}(qz)$$
transforms  $g(z)$  to a $z$ independent diagonal form, up to the action
of the complex affine Weyl group $\hat{W}$. Let $W$ be the $A_{N-1}$
Weyl group (the permutations of the Cartan elements). Then
 $\hat{W}= ({\bf Z}R^{\vee}\tau+{\bf Z}R^{\vee})\smd W$
($R^{\vee}$ is the dual root system).  The moduli space ${\cal L}$ 
in this case is 
the Weyl alcove. The comparison of two description  of
 holomorphic structures on elliptic curves (\ref{a25}) and (\ref{a12})
 was carried out
in \cite{EF,EKh} in terms of two loop current algebras and invariants of
their coadjoint actions.

In general case $(g>1)$ the gauge transform (\ref{a23}) allows to choose 
$g_a$ as constant matrices. They are defined up to the common conjugation by
a $\GL$ matrix. Thus the moduli space of holomorphic bundles in the 
(\ref{a24}) description are defined as the quotient
$$
{\cal L}\sim(\underbrace{\GL\oplus\ldots\oplus\GL}_{g})/\GL.
$$
Since the center of $\GL$  acts trivially we obtain  $\dim{\cal L}=
N^2(g-1)+1$ (see (\ref{a14})).

\section{Symplectic geometry in the double coset picture}
\setcounter{equation}{0}

Here we consider the Hitchin integrable systems which are defined on
the cotangent bundle $T^*{\cal L}$ 
to the moduli of stable holomorphic bundles ${\cal L}$. As it was done in the
original work
\cite{H} this space is derived as a symplectic quotient of 
$T^*{\cal L}_{\Si}^D$ under the gauge action of ${\cal G}_{\Si}^{C^{\infty}}$.
We will come to the same systems by the three step symplectic reductions from
some big upstairs space. The main object of this section is
 the commutative diagram (\ref{b8}), which describes these reductions and
intermediate spaces.

{\bf 1}. First, as intermediate step, consider the Hitchin description of 
$T^*{\cal L}$.
The upstairs phase space is the cotangent bundle $T^*{\cal L}_{\Si}^D$ to the
space ${\cal L}_{\Si}^D$ (\ref{a12}) of holomorphic structures
 on the bundle $P$
in the Dolbeault picture.
It is the space of pairs
\beq{b1}
T^*{\cal L}_{\Si}^D=
\{\phi,~d''_A,~~~\phi\in\Om_{C^{\infty}}^{(1,0)}(\Si,({\rm End}V)^*)\} .
\eq
The field $\phi$ is called the Higgs field and the bundle $T^*{\cal L}_{\Si}^D$ 
is the Higgs bundle. 
We can consider the Higgs field as a form
$$\phi\in\Om_{C^{\infty}}^{0}(\Si,({\rm End}V)^*\otimes K),$$
where $K$ is the canonical bundle on $\Si$.
Locally on ${\cal U}_a$ 
$$
d''_a=\bp+\bar A_a,~\bar A_a=h^{-1}_a\bp h_a ,
~h_a\in\Om_{C^{\infty}}^0({\cal U}_a,{\cal R}_\Si).
$$
The symplectic form on it is
\beq{b2}
\om^D=\int_{\Si}\tr(D\phi,Dd''_A).
\eq
The action of the gauge group ${\cal G}_{\Si}^{C^{\infty}}$ (\ref{a11}) on
$d_{ A}''$ (\ref{a12a}) with 
$$
\phi\rar\vf^{-1}\phi\vf
$$
is a symmetry of $T^*{\cal L}_{\Si}^D$. It defines the moment map
$$
\mu_{{\cal G}_{\Si}^{C^{\infty}}}(\phi,\bar A):~T^*{\cal L}_{\Si}^D\rar
{\rm Lie}^*({\cal G}_{\Si}^{C^{\infty}}),
$$
$$\mu_{{\cal G}_{\Si}^{C^{\infty}}}(\phi,\bar A)=[d_{A}'',\phi].
$$
For the zero level moment map $[d_{A}'',\phi]=0$  the  Higgs field becomes
holomorphic 
$$\phi\in H^0(\Si,({\rm End}V)^*\otimes K) .$$
The symplectic quotient
$\mu^{-1}(0)/{\cal G}_{\Si}^{C^{\infty}}=
T^*{\cal L}_{\Si}^D//{\cal G}_{\Si}^{C^{\infty}}$ 
is identified with the cotangent
bundle to the moduli space $T^*{\cal L}$. The Hitchin 
commuting integrals  are constructed
by means of $(1-j,1)$ holomorphic differentials $\nu_{j,k},k=1,\ldots$:
\beq{b3}
I^D_{j,k}=\int_{\Si}\nu^D_{j,k}\tr\phi^j.
\eq
Since the space of these differentials $H^0(\Si,K\otimes T^j)$ ($K$ is the
canonical class,\\ $T^j$ is $(-j,0)$ forms) 
has dimension $(2j-1)(g-1)$ for $j>1$ and $g$ for $j=1$ we have $N^2(g-1)+1$
independent commuting integrals, providing the complete integrability
of the Hamiltonian systems (\ref{b2}),(\ref{b3}). The integrals (\ref{b3}) 
define the Hitchin map
$$H^0(\Si,({\rm End}V)^*\otimes K)\rar H^0(\Si, K^j).$$
\bigskip

{\bf 2.} The same system can be derived starting from the cotangent bundle
$T^*{\cal L}_{\Si}^{Ch}$ to the holomorphic structures on $P$ defined 
in the C-type description (\ref{a16}). Now
\beq{b4}
T^*{\cal L}_{\Si}^{Ch}=\{\eta_{ab},g_{ab}|
~\eta_{ab}\in\Om_{hol}^{(1,0)}({\cal U}_{ab},({\rm End}V)^*),~
g_{ab}\in\Om_{hol}^0({\cal U}_{ab},P)\},
\eq
This bundle can be endowed with the  symplectic structure by means of 
the Cartan-Maurer one-forms on $\Om_{hol}^0({\cal U}_{ab},P)$. 
Let
$\G_a^b(C,D)$ be a path in ${\cal U}_{ab}$ with the end points
in the triple intersections 
$C\in {\cal U}_{abc}={\cal U}_{a}\cap{\cal U}_{b}\cap{\cal U}_{c}$,
$D\in {\cal U}_{abd}$. 
We can put the data (\ref{b4}) on the fat graph corresponding to the covering 
$\{{\cal U}_a\}$. 
The   edges of the graph are $\{\G_a^b(CD)\}$
and $\{\G_b^a(DC)\}$ with opposite orientation.
 We assume that the covering
is such that the orientation of edges defines the oriented contours
around the faces ${\cal U}_a$. The fields $\eta_{ab},g_{ab}$ are attributed
 to the edge $\G_a^b(CD)$, while $\eta_{ba},g_{ba}$ to $\G_b^a(DC)$.
 The last pair is not independent -  $(g_{ab}^{-1}=g_{ba})$  (see (\ref{a4})).
Its counterpart 
 in the dual space is
\beq{b4a}
\eta_{ab}(z_a)=g_{ab}(z_a)\eta_{ba}(z_b(z_a))g_{ab}^{-1}(z_a).
\eq
The symplectic structure is defined by the form
\beq{b5}
\om^{Ch}=\sum_{\rm edges}
\int_{\G_a^b(CD)}D\tr(\eta_{ab}(z_a)(Dg_{ab}g_{ab}^{-1})(z_a)).
\eq
Here the sum is taken over the edges of the oriented graph obtained
from the fat graph after the identifications of fields (\ref{b4a}).
 In other words
we consider only the edge $\G_a^b$ with the  fields $g_{ab},\eta_{ab}$
and forget about the edge $\G_b^a$.
Since $\eta_{ab}$ and $g_{ab}$ are holomorphic in ${\cal U}_{ab}$, 
the definition is
independent on a choice of the path $\G_a^b$ within ${\cal U}_{ab}$.
The symplectic form is invariant under the gauge transformations (\ref{a17})
 supplemented by
\beq{b6}
\eta_{ab}\rar f_a\eta_{ab}f_a^{-1}.
\eq
 The set of invariant 
commuting integrals on $T^*{\cal L}_{\Si}^{Ch}$ is
\beq{b7}
I^{Ch}_{j,k}=\sum_{\rm edges}\int_{\G_a^b(CD)}
\nu^{Ch}_{(j,k)}(z_a)\tr(\eta_{ab}^j(z_a)),
\eq
where $\nu^{Ch}_{j,k}$ are $(1-j,0)$ differentials, which are related locally 
to the $(1-j,1)$ differentials as $\nu^D_{j,k}=\bp \nu^{Ch}_{j,k}$.

We can consider the system on the defined above  graph $L_{ab}$
 which
is dual to $\G_a^b(CD)$. The fields $g_{ab},\eta_{ab}, a,b=1\ldots$ are
defined on edges, while the gauge transformations $f_a$ live on 
vertices.

The moment map is
$$
\mu_{{\cal G}_{\Si}^{hol}}(\eta_{ab},g_{ab}):~T^*{\cal L}_{\Si}^{Ch}\rar
{\rm Lie}^*({\cal G}_{\Si}^{hol}).
$$
According to (\ref{a21}) the Hamiltonian generating the gauge 
transformations is
$$
F_{\ep^{hol}}=
\sum_{\rm edges}
\int_{\G_a^b(CD)}
\tr(\eta_{ab}(z_a)\ep_a^{hol}(z_a))-
\tr(\eta_{ab}(z_a)g_{ab}(z_a)\ep_b^{hol}((z_b(z_a))g_{ab}(z_a)^{-1})=
$$
$$
\sum_{\rm edges}
\int_{\G_a^b(CD)}
\tr(\eta_{ab}(z_a)\ep_a^{hol}(z_a))-
\tr(\eta_{ba}(z_b(z_a))\ep_b^{hol}(z_b(z_a)))=
$$
$$
\sum_a\int_{\G_a}\sum_b\tr(\eta_{ab}(z_a)\ep_a^{hol}(z_a)),
$$
where $\G_a$ is an oriented contour around ${\cal U}_a$.
The moment equation $\mu_{{\cal G}_{\Si}^{hol}}=0$ can be read off
from $F_{\ep^{hol}}$. It means that $\eta_{ab}$
 is a boundary value of some
holomorpfic form defined on ${\cal U}_a$
\beq{b8a}
\eta_{ab}(z_a)=H_a(z_a),~{\rm for}~z_a\in{\cal U}_{ab},~H_a\in
\Om^{(1,0)}_{hol}({\cal U}_a,({\rm End}V)^*).
\eq
The reduced system is again the cotangent bundle to the moduli space of
holomorphic bundles
$$
T^*{\cal L}={\cal G}_\Si^{hol}\backslash\backslash T^*{\cal L}_{\Si}^{Ch}=
{\cal G}_\Si^{hol}\backslash\mu_{{\cal G}_{\Si}^{hol}}^{-1}(0),
$$
which has dimension $2N^2(g-1)+2$.
\bigskip

{\bf 3}. To get the cotangent bundle $T^*{\cal L}$ via the symplectic reduction 
we can start from  $T^*{\cal R}_\Si$ using the double coset 
representation (\ref{a15}). Then  $T^*{\cal L}_{\Si}^D$
or $T^*{\cal L}_{\Si}^{Ch}$ are obtained on the intermediate stages of the two step
reduction under the actions of ${\cal G}_\Si^{hol}$ or ${\cal G}_\Si^{C^\infty}$.
Since these groups act from different sides on ${\cal R}_\Si$ 
their actions commute and the result
of the reduction procedure is independent on their order. But the  space 
${\cal R}_\Si$, as we already have remarked, is not free - 
its elements satisfy (\ref{a6}).
We will represent the constraints (\ref{a6}) as moment constraints and
consider the "superfree" space - cotangent bundle to the group ${\cal G}_{\Si}$
(\ref{a8}). More exactly we will consider (Theorem 4.1) the three step 
symplectic reductions which result in
 the following commutative "tadpole" diagram 
\bigskip
\beq{b8}
\begin{array}{ccccc}
& &\fbox{$T^*{\cal G}_{\Si}$}& &\\
& {\cal G}_{\Si}^{A}&\downarrow& &\\
& &\fbox{$T^*{\cal R}_{\Si}$}& &\\
&{\cal G}_\Si^{hol}\swarrow& &\searrow{\cal G}_{\Si}^{C^{\infty}}&\\
\fbox{$T^*{\cal L}_{\Si}^D$}&
& & &\fbox{$T^*{\cal L}_\Si^{Ch}$}\\
&{\cal G}_{\Si}^{C^{\infty}}\searrow&&\swarrow{\cal G}_\Si^{hol}&\\
&&\fbox{$T^*{\cal L}$}&&
\end{array}
\eq
\bigskip

To begin with we define the initial data on $T^*{\cal G}_{\Si}$ and
the gauge group ${\cal G}_{\Si}^{A}$.
To construct $T^*{\cal G}_{\Si}$ we consider three dual elements
$$
\Psi_a\in\Om_{C^{\infty}}^{(1,1)}({\cal U}_{a},({\rm End}V)^*),~
\eta_{ab},\eta_{ba}\in\Om_{C^{\infty}}^{(1,0)}({\cal U}_{ab}, ({\rm End}V)^*),
$$
$$
\xi_{ab},\xi_{ba}\in\Om_{C^{\infty}}^{(0,1)}({\cal U}_{ab},({\rm End}V)^*).
$$
Cotangent bundle $T^*{\cal G}_{\Si}$ is the set of fields 
$(\Psi,\eta,\xi,h)$. We endow it with the symplectic structure. 
Consider the same fat graph with edges $\G_a^b(CD)$ and $\G_b^a(DC)$ as
in {\bf 4.2}.
Then
 \beq{b9}
\om_{\Si}=D\{\sum_a\int_{{\cal U}_a}\tr(\Psi_aDh_ah_a^{-1})+
\sum_{b}[\int_{\G_a^b}\tr(\eta_{ab}Dh_ah_a^{-1})+
\int_{\G_{a}^b}\tr(\xi_{ab}h_a^{-1}Dh_a)]\}.
\eq
We assume as before that
 paths $\G_a^b,\G_a^c,\dots$ can be unified in 
a closed oriented contour $\G_a\subset{\cal U}_a$ .
 The integral over ${\cal U}_{a}$ means 
in fact the integral over a part of ${\cal U}_{a}$  restricted by
the contour $\G_a$. Thus the first sum can be replaced by
the integration over $\Si$. To maintain the independence of $\om_{\Si}$
on the choice of the contours $\G_a$ we introduce the following "gauge"
symmetry. Its action defines of variations of fields along with  variations
of contours. Let $\G'_a$ be another contour and $\de{\cal U}_{a}$ be the
corresponding variation of the integration domain. There is the 
integral relation between fields coming from the Stokes theorem,
 providing the independence of $\om_{\Si}$
\beq{b9a}
\int_{\de{\cal U}_a}\tr(\Psi_aDh_ah_a^{-1})=
[\int_{\G'_a}\tr(\eta_{ab}Dh_ah_a^{-1})+
\int_{\G'_{a}}\tr(\xi_{ab}h_a^{-1}Dh_a)]\eq
$$
-[\int_{\G_a}\tr(\eta_{ab}Dh_ah_a^{-1})+
\int_{\G_a}\tr(\xi_{ab}h_a^{-1}Dh_a)].
$$
In other words, the variation of contours is compensated by the variation of
the field $\Psi$.

The form $\om_{\Si}$ (\ref{b9}) is invariant under the actions of 
${\cal G}_\Si^{hol}$ :
\beq{b12}
h_a\rar f_ah_a,~\Psi_a\rar f_a\Psi_af_a^{-1},~
\eta_{ab}\rar f_a\eta_{ab}f_a^{-1},~
\xi_{ab}\rar\xi_{ab},
\eq
and ${\cal G}_{\Si}^{C^{\infty}}$
\beq{b13}
h_a\rar h_a\vf_a,~\Psi_a\rar\Psi_a,~\eta_{ab}\rar\eta_{ab},~
\xi_{ab}\rar\vf_a^{-1}\xi_{ab}\vf_a.
\eq
We extend these group transformations by the following affine action of
the group
\beq{b10}
{\cal G}_{\Si}^{A}=
\{s_{ab}\in\Om_{C^{\infty}}^0({\cal U}_{ab},P)|s_{ab}=s_{ba}\}
\eq
 on $\xi_{ab}$
$$
\xi_{ab}\rar\xi_{ab}-s_{ab}^{-1}(\bp+h_a^{-1}\bp h_a)s_{ab}+h_a^{-1}\bp h_a
$$
leaving the other fields untouched. This action commutes with 
${\cal G}_{\Si}^{hol}$,
 but does not commute with ${\cal G}_{\Si}^{C^{\infty}}$. ${\cal G}_{\Si}^{A}$
 can be imbedded in  the bicomplex ${\bf G}$ (see (\ref{b10})).
On the Lie algebra level 
we have
\beq{b11}
\xi_{ab}\rar\xi_{ab}-(\bp\ep_{ab}^A +[h_a^{-1}\bp h_a,\ep^A_{ab}])
\eq
$$
(\ep_{ab}^A\in{\rm Lie}({\cal G}_{\Si}^{A})
=\{\Om_{C^{\infty}}^0({\cal U}_{ab},{\rm End}V)|~\ep_{ab}^A=\ep_{ba}^A\}.
$$ 
\bigskip

\begin{predl}
The form $\om_{\Si}$ (\ref{b9}) is invariant under the action of 
${\cal G}_{\Si}^{A}$.
\end{predl}
{\sl Proof}. From (\ref{b11}) 
$$\de_{\ep_{ab}^A}\xi_{ab}=-(\bp\ep_{ab}^A+[h_a^{-1}\bp h_a,\ep_{ab}^A]),$$
where $\ep^A\in{\rm Lie}({\cal G}_{\Si}^{A})$. Then
$$
-\de_{\ep^A}\om_{\ti{\Si}}=
\sum_{a}\sum_b
\int_{\G_a^b}\tr\{D([h_a^{-1}\bp h_a,\ep_{ab}^A])h_a^{-1}Dh_a
+\bp\ep_a^AD(h_a^{-1}Dh_a)+[h_a^{-1}\bp h_a,\ep_a^A]D(h_a^{-1}Dh_a)\}
$$
$$
=-\sum_{a}\sum_b
\int_{\G_a^b}\tr\{[D(h_a^{-1}Dh_a),h_a^{-1}Dh_a]+
\bp D(h_a^{-1}Dh_a)+[h_a^{-1}\bp h_a,D(h_a^{-1}Dh_a)],\ep_{ab}^A\}.
$$
Then direct calculations show that the sum  under the integral in front
of $\ep_{ab}^A$ vanishes.
 Therefore $\om_{\Si}$ is invariant under these
transformations. $\Box$
\bigskip

More generally,  the dual fields 
$(\Psi_a,\eta_{ab},\eta_{ba},\xi_{ab},\xi_{ba})$
 can be incorporated in a general pattern of two exact
${\bf G}$ bimoduli:
 
\begin{center}
${\bf M'}^*$
\end{center}
$$
\def\normalbaselines{\baselineskip20pt
       \lineskip3pt    \lineskiplimit3pt}
\def\mapright#1{\smash{
        \mathop{\longrightarrow}\limits^{#1}}}
\def\mapup#1{\Big\uparrow\rlap
       {$\vcenter{\hbox{$\scriptstyle#1$}}$}}
\matrix{
&\mapup{\p^C}&   &\mapup{\p^C}&   &\mapup{\p^C}&\cr
0\rar&\Om_{hol}^{(1,0)}({\cal U}_{ab}, ({\rm End}V)^*)&\mapright{i} 
&\Om_{C^{\infty}}^{(1,0)}({\cal U}_{ab},({\rm End}V)^*) & \mapright{\bp} & 
\Om_{C^{\infty}}^{(1,1)}({\cal U}_{ab},({\rm End}V)^*)& \mapright{\bp}\cr
&\mapup{\p^C}&   &\mapup{\p^C}&   &\mapup{\p^C}&\cr
0\rar&\Om_{hol}^{(1,0)}({\cal U}_{a},({\rm End}V)^*)
&\mapright{i} 
&\Om_{C^{\infty}}^{(1,0)}({\cal U}_{a},({\rm End}V)^*)& \mapright{\bp} & 
\Om_{C^{\infty}}^{(1,1)}({\cal U}_{a},({\rm End}V)^*)
& \mapright{\bp}\cr
&\mapup{i}&   &\mapup{i}&   &\mapup{i}&\cr
0\rar&\Om_{hol}^{(1,0)}(\Si,({\rm End}V)^*)&\mapright{i} 
&\Om_{C^{\infty}}^{(1,0)}(\Si,({\rm End}V)^*)& \mapright{\bp} & 
\Om_{C^{\infty}}^{(1,1)}(\Si,({\rm End}V))^*&
 \mapright{\bp}\cr
&\mapup{}&   &\mapup{}&   &\mapup{}&\cr
&0&            &0&             &0&\cr}
$$

\bigskip

\begin{center}
 ${\bf M''}^*$
\end{center}
$$
\def\normalbaselines{\baselineskip20pt
       \lineskip3pt    \lineskiplimit3pt}
\def\mapright#1{\smash{
        \mathop{\longrightarrow}\limits^{#1}}}
\def\mapup#1{\Big\uparrow\rlap
       {$\vcenter{\hbox{$\scriptstyle#1$}}$}}
\matrix{
&\mapup{\p^C}&   &\mapup{\p^C}&   &\mapup{\p^C}&\cr
0\rar&\Om_{antihol}^{(0,1)}({\cal U}_{ab}, ({\rm End}V)^*)&\mapright{i} 
&\Om_{C^{\infty}}^{(0,1)}({\cal U}_{ab},({\rm End}V)^*) & \mapright{\p} & 
\Om_{C^{\infty}}^{(1,1)}({\cal U}_{ab},({\rm End}V)^*)& \mapright{\p}\cr
&\mapup{\p^C}&   &\mapup{\p^C}&   &\mapup{\p^C}&\cr
0\rar&\Om_{antihol}^{(0,1)}({\cal U}_{a},({\rm End}V)^*)
&\mapright{i} 
&\Om_{C^{\infty}}^{(0,1)}({\cal U}_{a},({\rm End}V)^*)& \mapright{\p} & 
\Om_{C^{\infty}}^{(1,1)}({\cal U}_{a},({\rm End}V)^*)
& \mapright{\p}\cr
&\mapup{i}&   &\mapup{i}&   &\mapup{i}&\cr
0\rar&\Om_{antihol}^{(0,1)}(\Si,({\rm End}V)^*)&\mapright{i} 
&\Om_{C^{\infty}}^{(0,1)}(\Si,({\rm End}V)^*)& \mapright{\p} & 
\Om_{C^{\infty}}^{(1,1)}(\Si,({\rm End}V))^*&
 \mapright{\p}\cr
&\mapup{}&   &\mapup{}&   &\mapup{}&\cr
    &0&        &0&          &0&\cr}
$$ 

\bigskip
We remind that 
$$\Psi_a\in\Om_{C^{\infty}}^{(1,1)}({\cal U}_{a},({\rm End}V)^*),~
\eta_{ab},\eta_{ba}\in\Om_{C^{\infty}}^{(1,0)}({\cal U}_{ab},({\rm End}V)^*),
$$
$$ 
\xi_{ab},\xi_{ba}\in\Om_{C^{\infty}}^{(0,1)}({\cal U}_{ab},({\rm End}V)^*).
$$
We will see that after the symplectic reductions these fields will obey some
 special constraints. Now we have all initial data 
to start from the top of the diagram (\ref{b8}) -the fields, the symplectic form
$\om_\Si$ (\ref{b9}) and the gauge groups actions 
(\ref{b12}),(\ref{b13}),(\ref{b10}).
\bigskip
\begin{theor}
There exist two ways of symplectic reductions represented by the
commutative diagram (\ref{b8}) which leads from $T^*{\cal G}_{\Si}$
to the cotangent bundle to the moduli space $T^*{\cal L}$.
\end{theor}
To prove Theorem we shall go down along the diagram.
\bigskip

{\bf 4}.
Consider first the action of ${\cal G}_{\Si}^{A}$ (\ref{b11}).
Let $T^*{\cal R}_{\Si}=\{\Psi,\eta,h\}$ and $h$ is GTF
with symplectic form 
\beq{b14}
\om_{\Si}=D\{\sum_a[\int_{{\cal U}_a}\tr(\Psi_aDh_ah_a^{-1})+
\sum_{b}\int_{\G_a^b(CD)}\tr(\eta_{ab}Dh_ah_a^{-1})]\}.
\eq

\begin{lem}
$$T^*{\cal R}_{\Si}=T^*{\cal G}_{\Si}//{\cal G}_{\Si}^{A}=
\mu_{{\cal G}_{\Si}^{A}}^{-1}(0)/{\cal G}_{\Si}^{A}.
$$
\end{lem}
{\sl Proof}. It follows from (\ref{b9}),(\ref{b10}) that
the Hamiltonian of ${\cal G}_{\Si}^{A}$ action is 
$$
F_{\ep^A}=\sum_{a>b}F_{ab},~F_{ab}=
\int_{\G_a^b(CD)}\tr(\ep_{ab}^Ah_a^{-1}\bp h_a)+
\int_{\G_b^a(DC)}\tr(\ep_{ba}^Ah_b^{-1}\bp h_b).
$$
In fact the one-form 
$$
DF_{ab}=\int_{\G_a^b(CD)}\{\tr(\bp\ep_{ab}^Ah_a^{-1}D h_a)+
\tr(\ep_{ab}^A[h_a^{-1}\bp h_a,h_a^{-1}D h_a])\}
$$
can be obtain from $\om_\Si$ (\ref{b9}) by the action of the
vector field generated by ${\cal G}_{\Si}^{A}$ (\ref{b10}).
But 
$\ep_{ab}^A=\ep_{ba}^A$ (\ref{b10}). Putting the moment equal to zero
$\mu_{{\cal G}_{\Si}^{A}}=0$ we come to the constraints
$h_a^{-1}\bp h_a=h_b^{-1}\bp h_b$, which are exactly (\ref{a6}).
Note that the gauge transform (\ref{b11}) allows to fix $\xi_{ab}=0$.
Thus  the symplectic quotient 
$T^*{\cal R}_\Si=T^*{\cal G}_\Si//{\cal G}_{\Si}^{A}$
 has the field content 
$(\Psi,\eta,h\in\Om_{C^{\infty}}^0(\Si,{\cal P}))$ with $\om_{\Si}$
(\ref{b14}). $\Box$
\bigskip

{\bf 5.}
Consider the action of ${\cal G}_\Si^{hol}$ (\ref{a21}),(\ref{b12}) on 
$T^*{\cal R}_\Si$,
which corresponds to the left arrow in the diagram (\ref{b9}). We will prove
\bigskip
\begin{lem}
$$
T^*{\cal L}_{\Si}^D={\cal G}_\Si^{hol}\backslash\backslash T^*{\cal R}_\Si=
{\cal G}_\Si^{hol}\backslash\mu_{{\cal G}_\Si^{hol}}^{-1}(0),
$$
where $T^*{\cal L}_{\Si}^D$ is defined by (\ref{b1}) with the symplectic 
structure (\ref{b2}). 
\end{lem}
{\sl Proof}. From (\ref{b12}) and (\ref{b14}) we read off the hamiltonian
 of the gauge fields
$$
F_{\ep^{hol}}=
\sum_{a}[
\int_{{\cal U}_{a}}\tr(\Psi_a\ep_a^{hol})+
\sum_{b}
\int_{\G_a^b}\tr(\eta_{ab}\ep_a^{hol})].
$$
On ${\cal U}_a$ we can put  $\Psi_a=\bp(\ti{\Phi}_a+H_a)$, where 
$\ti{\Phi}_a\in\Om_{C^{\infty}}^{(1,0)}({\cal U}_a,({\rm End}V)^*$ and $H_a$
is an arbitrary element from
 $\Om_{hol}^{(1,0)}({\cal U}_{a},({\rm End}V)^*)$ (see 
${\bf M'}^*$). Then 
$$
F_{\ep^{hol}}=
\sum_{a}\sum_{b}
\int_{\G_a^b}\tr((\ti{\Phi}_a+H_a+
\eta_{ab})\ep_a^{hol}).
$$
Resolving the moment constraint $\mu_{{\cal G}_\Si^{hol}}=0$ gives
\beq{b14d}
\eta_{ab}(z_a,\bar{z}_a)=
-\ti{\Phi}_a(z_a,\bar{z}_a)-H_a(z_a),~z_a\in{\cal U}_{ab}.
\eq
By means of the Stokes theorem $\om_{\Si}$ (\ref{b14}) can be transformed  
to the form
$$
\om_{\Si}=
D\{
\sum_{a}[\int_{{\cal U}_{a}}\tr(\bp(\ti{\Phi}_a+H_a)Dh_ah_a^{-1})+
\sum_{b}\int_{\G_a^b}\tr(\eta_{ab}Dh_ah_a^{-1})]\}=
$$
$$
-\sum_{a}D\int_{{\cal U}_{a}}\tr((\ti{\Phi}_a+H_a)\bp (Dh_ah_a^{-1})).
$$
 
Let 
\beq{b14c}
\phi_a=-h^{-1}_a(\ti{\Phi}_a+H_a)h_a . 
\eq
Remind that $H_a=H_a(z_a)$ is an
arbitrary holomorphic function on ${\cal U}_a$. 
We will choose it in a such way that  $\phi_a$ becomes a global section in
$ \Om_{C^{\infty}}^{(1,0)}(\Si,({\rm End}V)^*) $. In other words
\beq{b14e}
(\phi_a-\phi_b)|_{\G_a^b}=0.
\eq
In fact, since $g_{ab}=h_ah_b^{-1}$, 
\beq{b14b}
h_a(\phi_b-\phi_a)h^{-1}_a=(\ti{\Phi}_a-g_{ab}\ti{\Phi}_bg_{ab}^{-1})
+(H_a-g_{ab}H_bg_{ab}^{-1}),
\eq
where the second term is holomorphic.
Consider the integral $I_a$ over the contour $\G_a=\cup_b\G_a^b$ around 
${\cal U}_a$
$$
I_a=-\int_{\G_a}\sum_b\frac
{(\ti{\Phi}_a-g_{ab}\ti{\Phi}_bg_{ab}^{-1})(y)}
{z-y}dy.
$$
Due to the Sokhotsky-Plejel theorem \cite{Mu} $I_a$ is holomorphic
inside and outside $\G_a$. It has a jump 
$\ti{\Phi}_a-g_{ab}\ti{\Phi}_bg_{ab}^{-1}$ on the contour. 
Let 
$$
H_a=I_a ~{\rm in}~{\cal U}_a,
$$
$$
H_b=g_{ab}^{-1}I_ag_{ab}~{\rm outside}~{\cal U}_a.
$$
Therefore the functions
$H_a$ and $g_{ab}H_bg_{ab}^{-1}$ defining by this integral provide the 
vanishing of the left hand side (\ref{b14b}). 
The symplectic form
$\om_{\Si}$ in terms of $\phi$ and $\bar{A}$ can be rewritten as 
$$
\om_{\Si}= \sum_{a}D\int_{{\cal U}_{a}}\tr(\phi_a D\bar{A}_a)=
D\int_{\Si}\tr(\phi D\bar{A}).
$$
This form coincides with $\om^D$ (\ref{b2}) for $T^*{\cal L}_\Si^D$. 
The field $\phi$ is invariant under
the ${\cal G}_\Si^{hol}$ action (\ref{b12}). 
Therefore the symplectic reduction by the gauging  ${\cal G}_\Si^{hol}$
leaves us with the fields $\phi$ and $h$ and the symplectic structure 
(\ref{b2}).
In other words $T^*{\cal L}_\Si//{\cal G}_\Si^{hol}=T^*{\cal L}_\Si^D$.
 $\Box$
\bigskip
We can now move down along the left side of  diagram (\ref{b9}) as it
was described in {\bf 1.} and obtain eventually $T^*{\cal L}$.

It will be instructive to look on relations between two type 
of dual fields
$\eta$ (\ref{b14d}) and $\phi$ (\ref{b14c}) that arise after
 these two consecutive reductions.
 On the first step we found that $\eta$ are
boundary valued forms
\beq{b15a}
\eta_{ab}(z_a,\bar{z}_a)=h_a^{-1}\phi h_a|_{{\cal U}_{ab}}.
\eq
Moreover, it follows from (\ref{b14e}) that
\beq{b15b}
\eta_{ab}(z_a,\bar{z}_a)=
g_{ab}((z_a,\bar{z}_a)\eta_{ba}(z_b(z_a),\bar{z}_b(\bar{z}_a))
g_{ab}((z_a,\bar{z}_a)^{-1}.
\eq
The second reduction gives $\bp\phi+[\bar{A},\phi]=0$ (see {\bf 1.}).
It is equivalent to $\bp\eta=0$, due to (\ref{b15a}). 

\bigskip

{\bf 6}.
Now look on the right side of the  diagram. 
\begin{lem}
$$
T^*{\cal L}_{\Si}^{Ch}=T^*{\cal R}_\Si//{\cal G}_{\Si}^{C^{\infty}}=
\mu_{{\cal G}_{\Si}^{C^{\infty}}}^{-1}(0)/{\cal G}_{\Si}^{C^{\infty}},
$$
where $T^*{\cal L}_{\Si}^{Ch}$ is the cotangent bundle (\ref{b4}) 
with $\om^{Ch}$
(\ref{b5}).
\end{lem}
{\sl Proof.} 
The gauge action of ${\cal G}_{\Si}^{C^{\infty}}$
  (\ref{b13}) on  $T^*{\cal R}_\Si$ defines 
the Hamiltonian (see (\ref{b14}))
\beq{b14a}
F_{\ep^{C^{\infty}}}=
\sum_a\{\int_{{\cal U}_{a}}\tr(h_a^{-1}\Psi_ah_a\ep^{C^{\infty}}_a)+
\sum_{b}\int_{\G_a^b}\tr(h_a^{-1}\eta_{ab}h_a\ep^{C^{\infty}}_a)\}.
\eq
Consider the zero level of the moment map
$$
\mu_{{\cal G}_{\Si}^{C^{\infty}}}:~
T^*\ti{\cal L}_\Si\rar {\rm Lie}^*({\cal G}_{\Si}^{C^{\infty}}).
$$
From the first  terms in  (\ref{b14a}) we obtain
$$
\Psi_a=0,~a=1,\ldots.
$$
This choice of $\Psi$ breaks the invariance with respect 
to replacements of contours . But if $\bp\eta_{ab}=0$ then the exact form
of the path $\G_a^b(C,D)$ is nonessential. Note that this choice 
is consistent with
the definition of $\eta$ (\ref{b14e}) $(\eta_{ab}=H_a$ in the 
${\cal G}_\Si^{hol}$ reduction).
Picking up in the second sum in (\ref{b14a}) integrals over two 
neighbor edges
 we come to the condition
$$
\int_{\G_a^b}\tr(h_a^{-1}\eta_{ab}h_a\ep^{C^{\infty}}_a)+
\int_{\G_b^a}\tr(h_b^{-1}\eta_{ba}h_b\ep^{C^{\infty}}_b)=0.
$$
Since $\ep^{C^{\infty}}\in{\rm Lie}({\cal G}_{\Si}^{C^{\infty}})$ it is
"periodic" $\ep^{C^{\infty}}_a(z_a)=\ep^{C^{\infty}}_b(z_b(z_a))$. 
It gives the following form of constraints
$$
(h_a^{-1}\eta_{ab}h_a)(z_a)=(h_b^{-1}\eta_{ba}h_b(z_b(z_a))),
~z_a\in{\cal U}_{ab},
$$
or
\beq{b15}
\eta_{ab}(z_a)=g_{ab}(z_a)\eta_{ba}(z_b(z_a))g^{-1}_{ab}(z_a), ~
(g_{ab}(z_a)=h_a(z_a,\bar{z}_a)h_b^{-1}(z_b(z_a),\overline{z_b(z_a)})),
\eq
which is just the twisting property of $\eta$  (\ref{b4a}). Furthermore, 
the symplectic form 
$\om_{\Si}$ (\ref{b14}) due to vanishing the field $\Psi$ now is
$$
\om_{\Si}=\sum_{a>b}D[\int_{\G_a^b(CD)}\tr(\eta_{ab}Dh_ah_a^{-1})+
\int_{\G_b^a(DC)}\tr(\eta_{ba}Dh_bh_b^{-1})].
$$
Taking into account that
$$
(Dg_{ab}g_{ab}^{-1})(z_a)=
Dh_a(z_a)h_a^{-1}(z_a)-h_a(z_a)(h_b^{-1}Dh_b)(z_b(z_a))h_a(z_a)^{-1}
$$
and the moment constraint (\ref{b15}) we can rewrite $\om_{\Si}$ as
$$
\om_{\Si}=\sum_{\rm edges}
\int_{\G_a^b(CD)}D\tr(\eta_{ab}(z_a)(Dg_{ab}g_{ab}^{-1})(z_a)).
$$
It is just
$\om^C$ (\ref{b5}) in the 
C-type description of holomorphic bundles. We have the same
 field content and the same symplectic structure as in $T^*{\cal L}_{\Si}^{Ch}.$
  Therefore
$T^*{\cal R}_\Si//{\cal G}_{\Si}^{C^{\infty}}=T^*{\cal L}_{\Si}^{Ch}.$ $\Box$

The last step on the right side of diagram was described in {\bf 2}.
Its completes the proof of Theorem.

\section {Schottky description of Hitchin systems} 
\setcounter{equation}{0}
{\bf 1}. Now consider the last step in the diagram (\ref{b8})
 in the Schottky parameterization. Since in this case we have only one topologically 
nontrivial cell
$\ti{\Si}$ 
the symplectic reduction is differ from the described in {\bf 4.2} for
the standard covering.
In this case the holomorphic fields 
$\eta_a,g_a=g_a(z_a),~a=1,\ldots,g$ live in vicinities 
${\cal V}_a$ of
${\cal A}_a$-cycles, and $z_a$ are local parameters in the internal disks.
(see (\ref{a22})).
The phase space is
$$
T^*{\cal L}_{\Si}^{Ch}=\{\eta_{a},g_{a}|
~\eta_{a}\in\Om_{hol}^{(1,0)}({\cal V}_{a},({\rm End}V)^*),~
g_{a}\in\Om_{hol}^0({\cal V}_{a},P)\}.
$$
In other words in accordance with (\ref{a24})
$$
T^*{\cal L}_{\Si}^{Ch}=\oplus_{a=1}^g T^*L_a(GL(N,{\bf C})),
$$
and the loop groups $L_a(GL(N,{\bf C}))$ are extended by the
projective transformations of $z_a$ as in (\ref{a24}). 
The symplectic form on this object is (see(\ref{b5})
\beq{c1}
\om^{Ch}=\sum_{a=1}^gD\int_{{\cal A}_{a}}
\tr(\eta_{a}(z_a),Dg_{a}g^{-1}_{a}(z_a)).
\eq
The gauge transformations (\ref{a17}),(\ref{b6}) act as the common conjugations
 by global holomorphic
in $\ti{\Si}$ matrix functions
\beq{c2} 
\eta_a(z_a)\rar f(z_a)\eta_a f^{-1}(z_a),~~
g_a\rar f(z_a)g_a(z_a)f^{-1}(\ga_a(z_a)).
\eq
The invariant commuting Hamiltonians (\ref{b7}) in this parameterization
are 
\beq{c3}
I^C_{j,k}=\sum_{a}\int_{{\cal A}_{a}}\nu^C_{(j,k)}(z_a)\tr(\eta_{a}^j(z_a)),
\eq

The gauge transform (\ref{c2}) produces the moment map 
$\mu_{{\cal G}_{\Si}^{hol}}$,
which takes the form 
$$
\mu_{{\cal G}_{\Si}^{hol}}=
\eta_{a}(\ga_a(z_a))-(g_{a}^{-1}\eta_{a}g_{a})(z_a),~a=1,\ldots g.
$$
Assume as above that $\mu_{{\cal G}_{\Si}^{hol}}=0$:
\beq{b17}
\eta_{a}(\ga_a(z_a))-(g_{a}^{-1}\eta_{a}g_{a})(z_a)=0,~a=1,\ldots g,
\eq
which is twisting property (\ref{b4a}) in the Schottky picture. 

\bigskip 

{\bf 2}.
The solutions of the moment equations are known in a few degenerate cases
\cite{N}. We will consider here as an example of the above construction 
holomorphic bundles over elliptic curves with a marked point. 

Define the elliptic curve as the quotient
$$\Si_{\tau}={\bf C}^*/q^{\bf Z},~q=\exp 2\pi i\tau.~$$
In this case
$$T^*{\cal L}_{\Si}^{Ch}\sim(\eta(z),g(z);p,s)$$
where $s\in\GL$ is a group element in the marked point $z=1$ and 
$p\in {\rm Lie}^*(\GL)$. In addition to (\ref{c2})
$$
p\rar f(z)pf^{-1}(1),~s\rar f(1)s.
$$
The one form $\eta(z)$ has a pole in the singular point $z=1$.
The symplectic form (\ref{c1}) on these objects is
\beq{c4}
\om^{Ch}=D\int_{{\cal A}}
\tr(\eta(z)Dgg^{-1}(z))+D\tr(pDss^{-1}).
\eq
The transition map  $g(z)$ can be diagonalized by (\ref{c2}):
$$
g(z)=\exp 2\pi i \vec{u}=\exp\{\di ~ 2\pi i(u_1,\ldots, u_N)\},
$$
where $u_j$ are $z$-independent.
We keep the same notation for the transformed $\eta(z)=
\Si_{n\in {\bf Z}}\eta_{j,k}^{(n)}z^n$.
The moment equation (\ref{c4}) takes the form 
$$
\eta(qz)-(g^{-1}\eta g)(z)=p\de(z),
$$
Rewrite it as
$$
q^n\eta_{j,k}^{(n)}-e^{2\pi i(x_k-x_j)}\eta_{j,k}^{(n)}=p_{j,k}^{(n)}.
$$
After the resolving the moment constraints we find
$$\eta_{j,j}(z)=w_j,~p_{jj}=0,$$
$$\eta_{j,k}=-\frac
{1}{2\pi i }
\frac
{\te(u_j-u_k-\ze)\te'(0)}
{\te(u_j-u_k)\te(\ze)},~z=\exp 2\pi i \ze,$$
where $w_j$ are new free parameters and
$\te(\ze)=\sum_{n\in {\bf Z}}e^{\pi i(n^2\tau+2n\ze)}$.
The symplectic form (\ref{c4}) on the reduced space takes the form
$$
\om^{red}=D\vec{w}\cdot D\vec{u} +
\tr D(J s^{-1}Ds),
$$
and $J$ defines the coadjoint orbit $p=s^{-1}Js$.
Consider the quadratic Hamiltonian (\ref{c3}).
 After the reduction $H$ takes
the form of the N-body elliptic Calogero Hamiltonian
with the spins \cite{Kr}:
$$
H= \frac{1}{2}(\vec{w}\cdot \vec{w}+\frac{1}{4\pi^2}
\sum_{j>k}^N[p_{j,k}p_{k,j}\wp(u_j-u_k|\tau)+E_2(\tau)]).
$$
Here $E_2(\tau)$ is the normalized Eisenstein series.

\bf Acknowledgments. 
{\sl We would like to thank  V.Fock,  B.Khesin, N.Nekrasov and A.Rosly
   for illuminating discussions.
 We are grateful to  the Max Planck Institute for Mathematik in Bonn 
for the hospitality where this work was prepared. The work of A.L.
is supported in part by the grant INTAS 944720 and the grant ISF NSR-300.
The work of M.O. is supported in part by the grant 
CEE-INTAS 932494 and the grant RFFI-96-02-18046 }

\small{

}

\end{document}